\documentclass{amsart}
\usepackage[dvips]{graphicx}
\usepackage{ifthen,amsfonts,amsmath,amssymb,epic,eepic}  
\usepackage{epsfig,color}  

\theoremstyle{plain}
\newtheorem{theorem}{Theorem}[section]
\newtheorem{lemma}[theorem]{Lemma}
\newtheorem{corollary}[theorem]{Corollary}
\newtheorem{proposition}[theorem]{Proposition}

\theoremstyle{remark}

\theoremstyle{definition}

\newcommand{\Energy}{{\text{Energy}}}

\newcommand{\Graph}{\operatorname{Graph}}

\def\ZZ{{\bold Z}}
\def\RP{{\bold RP}}
\def\RR{{\bold R}}
\def\SS{{\bold S}}

\newcommand{\nn}{{\bold {n}}}
\newcommand{\HH}{{\bold {H}}}

\newcommand{\dv}{{\text{div }}}

\newcommand{\Vol}{{\text {Vol}}}
\newcommand{\Area}{{\text {Area}}}
\newcommand{\Ric}{{\text {Ric}}}
\newcommand{\Scal}{{\text {Scal}}}

\def\CC{{\bold C }}
\newcommand{\eqr}[1]{(\ref{#1})}

\newcommand{\e}{{\text {e}}}
\newcommand{\cB}{{\mathcal {B}}}

\newcommand{\cS}{{\mathcal {S}}}

\newcommand{\cM}{{\mathcal {M}}}

\newcommand{\Conv}{\operatorname{Conv}}

\numberwithin{equation}{section}

\begin{document}

\title[Minimal submanifolds]
{Minimal submanifolds}

\author{Tobias H. Colding}%
\address{Department of Mathematics\\
MIT\\ 77 Mass. Ave.\\ Cambridge, MA 02139}
\author{William P. Minicozzi II}%
\address{Department of Mathematics\\
Johns Hopkins University\\
3400 N. Charles St.\\
Baltimore, MD 21218}
\thanks{The authors were partially supported by NSF Grants DMS
0104453 and DMS  0405695}

\email{colding@cims.nyu.edu and minicozz@math.jhu.edu}






\maketitle \tableofcontents

\section{Introduction}

Soap films, soap bubbles, and surface tension were extensively
studied by the Belgian physicist and inventor (the inventor of the
stroboscope) Joseph Plateau in the first half of the nineteenth
century.  At least since his studies, it has been known that the
right mathematical model for soap films are minimal surfaces --
the soap film is in a state of minimum energy when it is covering
the least possible amount of area.   Minimal surfaces and
equations like the minimal surface equation have served as
mathematical models for many physical problems.

The field of minimal surfaces dates back to the publication in
1762 of Lagrange's famous memoir ``Essai d'une nouvelle m\'ethode
pour d\'eterminer les maxima et les minima des formules
int\'egrales ind\'efinies''.  Euler had already in a paper
published in 1744 discussed minimizing properties of the surface
now known as the catenoid, but he only considered variations
within a certain class of surfaces. In the almost one quarter of a
millennium that has past since Lagrange's memoir minimal surfaces
has remained a vibrant area of research and there are many reasons
why.  The study of minimal surfaces was the birthplace of
regularity theory.  It lies on the intersection of nonlinear
elliptic PDE, geometry,   topology and general relativity.

In what follows we give a quick tour through many of the classical
results in the field of minimal submanifolds, starting at the
definition.

The field of minimal surfaces remain extremely active and has very
recently seen major developments that have solved many
longstanding open problems and conjectures; for more on this, see
the expanded version of this survey \cite{CM3}.  See also the
recent surveys \cite{MeP}, \cite{Pz} and the expository article
\cite{CM4}.

Throughout this survey, we refer to \cite{CM1} for references
unless otherwise noted.

\part{Classical and almost classical results}       \label{p:1}

Let $\Sigma \subset \RR^n$ be a \underline{smooth} $k$-dimensional
submanifold (possibly with boundary) and $C^{\infty}_0(N\Sigma)$
 the space of all infinitely differentiable,  compactly
supported, normal vector fields on $\Sigma$.  Given $\Phi$ in
$C^{\infty}_0(N\Sigma)$, consider the one--parameter variation
\begin{equation}
\Sigma_{t,\Phi}=\{x+t\,\Phi (x)  | x\in \Sigma\}\, .
\end{equation}
The so called first variation formula of volume is the equation
(integration is with respect to $d\text{vol}$)
\begin{equation}  \label{e:frstvar}
\left.\frac{d}{dt} \right|_{t=0}\Vol (\Sigma_{t,\Phi})
=\int_{\Sigma} \langle \Phi\, , \, \HH \rangle \, ,
\end{equation}
where $\HH$ is the mean curvature (vector) of $\Sigma$.  (When
$\Sigma$ is noncompact, then $\Sigma_{t,\Phi}$ in \eqr{e:frstvar}
is replaced by $\Gamma_{t,\Phi}$, where $\Gamma$ is any compact
set containing the support of $\Phi$.) The submanifold $\Sigma$ is
said to be a {\it minimal} submanifold (or just minimal) if
\begin{equation}
\left.\frac{d}{dt} \right|_{t=0}\Vol (\Sigma_{t,\Phi})=0
\,\,\,\,\,\,\,\,\,\,\,\text{ for all } \Phi\in
C^{\infty}_0(N\Sigma)
\end{equation}
 or, equivalently by \eqr{e:frstvar}, if the
mean curvature $\HH$ is identically zero.  Thus $\Sigma$ is
minimal if and only if it is a critical point for the volume
functional. (Since a critical point is not necessarily a minimum
the term ``minimal'' is misleading, but it is time honored.  The
equation for a critical point is also sometimes called the
Euler--Lagrange equation.)

Suppose now, for simplicity, that $\Sigma$ is an oriented
hypersurface with unit normal $\nn_{\Sigma}$.  We can then write a
normal vector field $\Phi \in C^{\infty}_0(N\Sigma)$  as $\Phi =
\phi \nn_{\Sigma}$, where function $\phi$ is in the space
$C^{\infty}_0(\Sigma)$ of infinitely differentiable, compactly
supported functions on $\Sigma$. Using this, a computation shows
that if $\Sigma$ is minimal, then
\begin{equation}    \label{e:secvar}
\left. \frac{d^2}{dt^2} \right|_{t=0}\Vol (\Sigma_{t,\phi
\nn_{\Sigma} }) =-\int_{\Sigma}\phi\,L_{\Sigma}\phi\, ,
\end{equation}
where
\begin{equation}
    L_{\Sigma}\phi=\Delta_{\Sigma}\phi+|A|^2\phi
\end{equation}
is the second variational (or Jacobi) operator. Here
$\Delta_{\Sigma}$ is the Laplacian on $\Sigma$ and $A$ is the
second fundamental form.  So $|A|^2=\kappa_1^2+\kappa_2^2+ \dots +
\kappa_{n-1}^2 $,
 where $\kappa_1,\dots \kappa_{n-1}$ are the principal curvatures of $\Sigma$
and $\HH= (\kappa_1+ \dots + \kappa_{n-1}) \, \nn_{\Sigma} $. A
minimal submanifold $\Sigma$ is said to be stable if
\begin{equation}
\left. \frac{d^2}{dt^2} \right|_{t=0}\Vol (\Sigma_{t,\Phi})\geq 0
\,\,\,\,\,\,\,\,\,\,\,\text{ for all } \Phi\in
C^{\infty}_0(N\Sigma)\, .
\end{equation}
Integrating by parts in \eqr{e:secvar}, we see that stability is
equivalent to the so called stability inequality
\begin{equation}   \label{e:stabin}
      \int |A|^2 \, \phi^2 \leq   \int |\nabla \phi|^2 \, .
\end{equation}
More generally, the {\it Morse index} of a minimal submanifold is
defined to be the number of negative eigenvalues of the operator
$L$.  Thus, a stable submanifold has Morse index zero.

\subsection{The Gauss map}

Let $\Sigma^2\subset \RR^3$ be a surface (not necessarily
minimal). The {\it Gauss map} is a continuous choice of a unit
normal $\nn:\Sigma\to \SS^2\subset \RR^3$. Observe that there are
two choices of such a map $\nn$ and $-\nn$ corresponding to a
choice of orientation of $\Sigma$. If $\Sigma$ is minimal, then
the Gauss map is an (anti) conformal map since the eigenvalues of
the Weingarten map are $\kappa_1$ and $\kappa_2=-\kappa_1$.
Moreover, for a minimal surface
\begin{equation}        \label{e:gmconf}
        |A|^2=\kappa_1^2+\kappa_2^2=-2\,\kappa_1\,\kappa_2=-2\,K_{\Sigma} \, ,
\end{equation}
where $K_{\Sigma}$ is the Gauss curvature.  It follows that the
area of the Gauss map is a multiple of the total curvature.

\subsection{Minimal graphs}

Suppose that $u:\Omega\subset \RR^2 \to \RR$ is a $C^2$ function.
The graph of   $u$
\begin{equation}
        \Graph_u= \{(x,y,u(x,y)) \mid (x,y)\in \Omega\}\, .
\end{equation}
has area
\begin{align}
        \Area (\Graph_u)&=\int_{\Omega} |(1,0,u_x)\times (0,1,u_y)|\\
        &=\int_{\Omega} \sqrt{1+u_x^2+u_y^2}=\int_{\Omega}\sqrt{1+|\nabla u|^2}
        \, ,\notag
\end{align}
and the (upward pointing) unit normal is
\begin{equation}                \label{e:bb1}
        \nn=\frac{(1,0,u_x)\times (0,1,u_y)}{|(1,0,u_x)\times (0,1,u_y)|}
        =\frac{(-u_x,-u_y,1)}{\sqrt{1+|\nabla u|^2}}\, .
\end{equation}
Therefore for the graphs $\Graph_{u+t\eta}$ where $\eta|\partial
\Omega=0$ we get that
\begin{equation}
        \Area (\Graph_{u+t\eta})
        =\int_{\Omega}\sqrt{1+|\nabla u+t\,\nabla \eta|^2}
\end{equation}
hence
\begin{align}
        \frac{d}{dt}_{t=0} \Area (\Graph_{u+t\eta})
        &=\int_{\Omega} \frac{\langle \nabla u \, , \nabla \eta \rangle}{\sqrt{1+|\nabla u|^2}}\\
&=-\int_{\Omega}\eta\, \dv\, \left( \frac{\nabla
u}{\sqrt{1+|\nabla u|^2}} \right)\, .\notag
\end{align}
It follows that the graph of $u$ is a critical point for the area
functional if and only if $u$ satisfies the divergence form
equation
\begin{equation}  \label{e:mineq}
        \dv \left( \frac{\nabla u}{\sqrt{1+|\nabla u|^2}} \right)=0  \, .
\end{equation}

Next we want to show that the graph of a function on $\Omega$
satisfying the minimal surface equation, i.e., satisfying
\eqr{e:mineq}, is not just a critical point for the area
functional but is actually area-minimizing amongst surfaces in the
cylinder $\Omega \times \RR \subset \RR^3$.  To show this, extend
first the unit normal $\nn$ of the graph in \eqr{e:bb1} to a
vector field, still denoted by $\nn$, on the entire cylinder
$\Omega \times \RR$. Let $\omega$ be the two-form on $\Omega
\times \RR$ given by that for $X,\, Y\in \RR^3$
\begin{equation}
\omega (X,Y)=\text{det} (X,Y,\nn)\, .
\end{equation}
An easy calculation shows that
\begin{equation}
d \omega=\frac{\partial }{\partial x} \left(
\frac{-u_x}{\sqrt{1+|\nabla u|^2}} \right) +\frac{\partial
}{\partial y} \left( \frac{-u_y}{\sqrt{1+|\nabla u|^2}} \right)
=0\, ,
\end{equation}
since $u$ satisfies the minimal surface equation. In sum, the
form $\omega$ is closed and, given any
  $X$ and $Y$ at a point $(x,y,z)$,
\begin{equation}        \label{e:calibration1}
        | \omega (X,Y) | \leq | X \times Y| \, ,
\end{equation}
where  equality holds if and only if
\begin{equation}        \label{e:calibration2}
        X, Y \subset T_{(x,y, u(x,y))}
        \Graph_u \, .
\end{equation}
  Such a form $\omega$ is called a {\it calibration}.
  From this, we have that if $\Sigma \subset \Omega \times \RR$
is any other surface with $\partial \Sigma=\partial\, \Graph_u$,
then by Stokes' theorem since $\omega$ is closed,
\begin{equation}  \label{e:areamin}
        \Area (\Graph_u) = \int_{\Graph_u} \omega =
                \int_{\Sigma} \omega \leq
        \Area (\Sigma) \,  .
\end{equation}
This shows that $\Graph_u$ is area-minimizing among all surfaces
in the cylinder and with the same boundary. If the domain $\Omega$
is convex, the minimal graph is absolutely area-minimizing. To see
this, observe first that if $\Omega$ is convex, then so is $\Omega
\times \RR $ and hence the nearest point projection $P: \RR^3 \to
\Omega \times \RR$ is a distance nonincreasing Lipschitz map that
is equal to the identity on
 $\Omega \times \RR$.
If $\Sigma \subset \RR^3$ is any other surface with $\partial
\Sigma=\partial\, \Graph_u$, then $\Sigma' = P(\Sigma)$ has $\Area
(\Sigma') \leq  \Area (\Sigma)$.  Applying \eqr{e:areamin} to
$\Sigma'$, we see that $\Area   (\Graph_u) \leq \Area (\Sigma')$
and the claim follows.

If $\Omega\subset \RR^2$ contains a ball of radius $r$, then,
since $\partial B_r\cap \Graph_u$ divides $\partial B_r$ into two
components at least one of which has area at most equal to $ (
\Area (\SS^2) / 2 ) \,r^2$, we get from \eqr{e:areamin} the crude
estimate
\begin{equation}                \label{e:area}
        \Area (B_r\cap \Graph_u)\leq \frac{\Area (\SS^2)}{2}\, r^2\, .
\end{equation}
When the domain $\Omega$ is convex, it is not hard to see that the
minimal graph is absolutely area-minimizing.

Very similar calculations to the ones above show that if
$\Omega\subset \RR^{n-1}$ and $u:\Omega\to \RR$ is a $C^2$
function, then the graph of $u$ is a critical point for the area
functional if and only if $u$ satisfies \eqr{e:mineq}. Moreover,
as in \eqr{e:areamin}, the graph of $u$ is actually
area-minimizing. Consequently, as in \eqr{e:area}, if $\Omega$
contains a ball of radius $r$, then
\begin{equation}                \label{e:areaforn}
        \Vol (B_r\cap \Graph_u)\leq \frac{\Vol (\SS^{n-1})}{2}\, r^{n-1}\, .
\end{equation}

\subsection{The maximum principle}

The first variation formula, \eqr{e:frstvar}, showed that a smooth
submanifold is a critical point for area if and only if the mean
curvature vanishes.  We will next
  derive the weak form of the  first variation formula
which is the basic tool for working with ``weak solutions''
(typically, stationary varifolds). Let $X$ be a vector field on
$\RR^n$.  We can write the divergence $\dv_{\Sigma} \, X$ of $X$
on $\Sigma$ as
\begin{equation}  \label{e:o1.3.1}
        \dv_{\Sigma} \, X=\dv_{\Sigma}\, X^T + \dv_{\Sigma}\, X^N = \dv_{\Sigma}\, X^T
        + \langle X , \HH \rangle \, ,
\end{equation}
where $X^T$ and $X^N$ are the tangential and normal projections of
$X$.  In particular, we get that, for a minimal submanifold,
\begin{equation}  \label{e:o1.3.1a}
        \dv_{\Sigma} \, X=\dv_{\Sigma}\, X^T
       \, .
\end{equation}
Moreover, from \eqr{e:o1.3.1} and Stokes' theorem, we see that
$\Sigma$ is minimal if and only if for all vector fields $X$ with
compact support and vanishing on the boundary of $\Sigma$,
\begin{equation}   \label{e:o1.3.2}
        \int_{\Sigma}\dv_{\Sigma}\, X=0\, .
\end{equation}
The key point is that \eqr{e:o1.3.2} makes sense as long as we can
define the divergence on $\Sigma$.   As a consequence of
\eqr{e:o1.3.2},
 we will show the following proposition:

\begin{proposition}
\label{p:haco} $\Sigma^k\subset \RR^n$ is minimal if and only if
the restrictions of the coordinate functions of $\RR^n$ to
$\Sigma$ are harmonic functions.
\end{proposition}

\begin{proof}
Let $\eta$ be a smooth function on $\Sigma$ with compact support
and $\eta |\partial \Sigma=0$, then
\begin{equation}
        \int_{\Sigma} \langle \nabla_{\Sigma} \eta , \,
        \nabla_{\Sigma} x_i \rangle = \int_{\Sigma} \langle  \nabla_{\Sigma}\eta \, ,
        e_i \rangle =  \int_{\Sigma} \dv_{\Sigma} (\eta\, e_i)\, .
\end{equation}
From this, the claim follows easily.
\end{proof}

Recall that if $\Xi\subset \RR^n$ is a compact subset, then the
smallest convex set containing $\Xi$ (the convex hull, $\Conv
(\Xi)$) is the intersection of all half--spaces containing $\Xi$.
The maximum principle forces a compact minimal submanifold to lie
in the convex hull of its boundary  (this is the  ``convex hull
property''):

\begin{proposition}
If $\Sigma^k\subset \RR^n$ is a compact minimal submanifold,
 then
$\Sigma\subset \Conv (\partial \Sigma)$.
\end{proposition}

\begin{proof}
A half--space $H \subset \RR^n$ can be written as
\begin{equation}
H = \{ x \in \RR^n \, | \, \langle x , e \rangle \leq a \} \, ,
\end{equation}
 for a vector $e \in \SS^{n-1}$ and constant $a \in \RR$.
  By Proposition \ref{p:haco},   the function
$u(x)=\langle e, x\rangle$ is harmonic on $\Sigma$ and hence
attains its maximum on $\partial \Sigma$ by the maximum principle.
\end{proof}

Another application of \eqr{e:o1.3.1a}, with a different choice of
vector field $X$, gives that for a $k$-dimensional minimal
submanifold $\Sigma$
\begin{equation}
        \Delta_{\Sigma} |x-x_0|^2=2 \, \dv_{\Sigma} (x-x_0)=2k \, .
\end{equation}
Later we will see that this formula plays a crucial role in the
monotonicity formula for minimal submanifolds.

\vskip2mm The argument in the proof of the convex hull property
can be rephrased as saying that as we translate a hyperplane
towards a minimal surface, the first point of contact must be on
the boundary.  When $\Sigma$ is a hypersurface, this is a special
case of the strong maximum principle for minimal surfaces:

\begin{lemma}             \label{l:smp}
Let $\Omega \subset \RR^{n-1}$ be an open connected neighborhood
of the origin. If $u_1$, $u_2:\Omega\to \RR$ are solutions of the
minimal surface equation with $u_1\leq u_2$ and $u_1(0)=u_2(0)$,
then $u_1\equiv u_2$.
\end{lemma}

Since any smooth hypersurface is locally a graph over a
hyperplane, Lemma \ref{l:smp} gives a maximum principle for smooth
minimal hypersurfaces.

 \vskip2mm  Thus far, the examples of minimal submanifolds have all been smooth.
The simplest non-smooth example
 is given by a pair of planes intersecting transversely along a
 line.  To get an example that is not even immersed, one can take
 three half--planes meeting along a line with an
 angle of
 $2\pi/3$ between each adjacent pair.

\section{Monotonicity and the mean value inequality}
\label{s:4}

Monotonicity formulas and mean value inequalities play a
fundamental role in many areas of geometric analysis.

\begin{proposition}
\label{p:pmonot} Suppose that $\Sigma^k\subset \RR^n$ is a minimal
submanifold and $x_0\in \RR^n$; then for all $0<s<t$
\begin{equation}        \label{e:monot0}
        t^{-k}\, \Vol (B_t (x_0) \cap \Sigma) -s^{-k}\, \Vol (B_s (x_0) \cap \Sigma)
        =\int_{(B_t(x_0) \setminus B_s(x_0) )\cap \Sigma} \frac{|(x-x_0)^N|^2}{|x-x_0|^{k+2}}\, .
\end{equation}
\end{proposition}

Notice that $(x-x_0)^N$ vanishes precisely when $\Sigma$ is
conical about $x_0$, i.e., when $\Sigma$ is invariant under
dilations about $x_0$. As a corollary, we get the following:

\begin{corollary}             \label{c:cmon}
Suppose that $\Sigma^k\subset \RR^n$ is a minimal submanifold and
$x_0\in \RR^n$; then the function
\begin{equation}                \label{e:thetadef}
        \Theta_{x_0}(s)
        =\frac{\Vol (B_s (x_0) \cap \Sigma)}{\Vol (B_s \subset \RR^k) }\,
\end{equation}
is a nondecreasing function of $s$.  Moreover, $\Theta_{x_0}(s)$
is constant in $s$ if and only if $\Sigma$ is conical about $x_0$.
\end{corollary}

Of course, if $x_0$ is a smooth point of $\Sigma$, then
$\lim_{s\to 0} \Theta_{x_0}(s) = 1$.  We will later see that the
converse is also true; this will be a consequence of the Allard
regularity theorem.

The monotonicity of area is a very useful tool in the regularity
theory for minimal surfaces --- at least when there is some {\it a
priori} area bound.  For instance, this monotonicity and a
compactness argument allow one to reduce many regularity questions
to questions about minimal cones (this was a key observation of W.
Fleming in his work on the Bernstein problem; see Section
\ref{s:5}).

\vskip2mm  Arguing as in Proposition \ref{p:pmonot}, we get a
weighted monotonicity:

\begin{proposition}
\label{p:meanvalue} If\hspace{2pt}  $\Sigma^k\subset \RR^n$ is a
minimal submanifold, $x_0\in \RR^n$, and $f$ is a function on
$\Sigma$, then
\begin{equation}
        t^{-k}\int_{B_t(x_0) \cap \Sigma} f-s^{-k}\int_{B_s(x_0) \cap \Sigma} f
\end{equation}
\begin{equation}
        =\int_{(B_t(x_0)\setminus B_s(x_0))\cap \Sigma} f\,
                \frac{|(x-x_0)^N|^2}{|x-x_0|^{k+2}}
        +\frac{1}{2} \int_{s}^t \tau^{-k-1}
        \int_{B_{\tau}(x_0)\cap \Sigma} (\tau^2-|x-x_0|^2)\,
                \Delta_{\Sigma} f\,d\tau
        \, . \notag
\end{equation}
\end{proposition}

We get immediately the following mean value inequality for the
special case of non--negative subharmonic functions:

\begin{corollary}
Suppose that $\Sigma^k\subset \RR^n$ is a minimal submanifold,
$x_0\in \RR^n$, and $f$ is a non--negative subharmonic function on
$\Sigma$;
 then
\begin{equation}                \label{e:defavg}
        s^{-k}  \int_{B_s(x_0) \cap \Sigma} f
\end{equation}
is a nondecreasing function of $s$. In particular, if $x_0\in
\Sigma$, then for all $s>0$
\begin{equation}
        f(x_0)\leq
        \frac{\int_{B_s (x_0) \cap \Sigma} f}
                {\Vol \, (B_s \subset \RR^k)} \, .
\end{equation}
\end{corollary}

\section{Rado's theorem}

One of the most basic questions is what does the boundary
$\partial \Sigma$ tell us about a compact minimal submanifold
$\Sigma$?  We have already seen that $\Sigma$ must lie in the
convex hull of $\partial \Sigma$, but there are many other
theorems of this nature.  One of the first is a beautiful result
of Rado which says that if $\partial \Sigma$ is a graph over the
boundary of a convex set in $\RR^2$, then $\Sigma$ is also graph
(and hence embedded). The proof of this uses basic properties of
nodal lines for harmonic functions.

\begin{theorem}            \label{t:rado}
Suppose that $\Omega\subset \RR^2$ is a convex subset and
$\sigma\subset \RR^3$ is a simple closed curve which is graphical
over $\partial \Omega$.  Then any minimal disk $\Sigma \subset
\RR^3$ with $\partial \Sigma =\sigma$ must be graphical over
$\Omega$ and hence unique by the maximum principle.
\end{theorem}

\begin{proof}
(Sketch.) The proof is by contradiction, so suppose that $\Sigma$
is such a minimal disk  and $x \in \Sigma$ is a point where
  the tangent plane to
$\Sigma$  is vertical.  Consequently, there exists $(a,b)\ne
(0,0)$ such that
\begin{equation}        \label{e:rado1}
        \nabla_{\Sigma} (a\,x_1+b\,x_2)(x) = 0 \, .
\end{equation}
By Proposition \ref{p:haco}, $a\,x_1+b\,x_2$ is harmonic on
$\Sigma$ (since it is a linear combination of coordinate
functions).  The local structure of nodal sets of harmonic
functions (see, e.g., \cite{CM1}) then gives that the level set
\begin{equation}                \label{e:rado2}
        \{ y \in \Sigma \, | \,  a\,x_1+b\,x_2 (y)
        = a\,x_1+b\,x_2 (x) \}
\end{equation}
 has a singularity
at $x$ where at least four different curves meet. If two of these
nodal curves were to meet again, then there would be a closed
nodal curve which must bound a disk (since $\Sigma$ is a disk).
By the maximum principle, $a\,x_1 + b\, x_2$ would   have to
 be constant on this disk
and hence constant on $\Sigma$ by unique continuation. This would
imply that $\sigma = \partial \Sigma$ is contained in the plane
given by \eqr{e:rado2}.  Since this is impossible, we conclude
that all of these curves go to the boundary without intersecting
again.

In other words, the plane in $\RR^3$ given by \eqr{e:rado2}
intersects $\sigma$ in at least four points.  However, since
 $\Omega \subset \RR^2$ is convex, $\partial \Omega$ intersects
 the line given
by \eqr{e:rado2} in exactly two points.  Finally, since $\sigma$
is graphical  over $\partial \Omega$, $\sigma$ intersects the
plane in $\RR^3$ given by \eqr{e:rado2} in exactly two points,
which gives the desired contradiction.
\end{proof}

\section{The theorems of Bernstein and Bers}        \label{s:5}

 A classical theorem of S. Bernstein from 1916 says that
entire (i.e., defined over all of $\RR^2$) minimal graphs are
planes. This remarkable theorem of Bernstein was one of the first
illustrations of the fact that the solutions to a nonlinear PDE,
like the minimal surface equation, can behave quite differently
from solutions to a linear equation.

\begin{theorem}   \label{t:bern} If
$u:\RR^{2}\to \RR$ is an entire solution to the minimal surface
equation, then $u$ is an affine function.
\end{theorem}

\begin{proof}
(Sketch.)  We will show that the curvature of the graph vanishes
identically; this implies that the unit normal is constant and,
hence, the graph must be a plane.  The proof follows by combining
two facts.  First, the area estimate for graphs \eqr{e:area} gives
\begin{equation}
        \Area (B_r\cap \Graph_u)\leq 2 \, \pi \, r^2\, .
\end{equation}
This quadratic area growth allows one to construct a sequence of
non-negative logarithmic cutoff functions $\phi_j$ defined on the
graph with $\phi_j \to 1$ everywhere and
\begin{equation}    \label{e:tozero}
   \lim_{j\to \infty} \,  \int_{\Graph_u} |\nabla \phi_j|^2 = 0 \, .
\end{equation}
Moreover, since graphs are area-minimizing, they must be stable.
We can therefore use $\phi_j$ in the stability inequality
\eqr{e:stabin} to get
\begin{equation}    \label{e:stabiii}
    \int_{\Graph_u} \phi_j^2 \, |A|^2 \leq \int_{\Graph_u} |\nabla \phi_j|^2  \, .
\end{equation}
Combining these gives that $|A|^2$ is zero, as desired.
\end{proof}

 Rather surprisingly,
this result very much depended on the dimension.  The combined
efforts of E. De Giorgi, F. J. Almgren, Jr., and J. Simons
 finally gave:

\begin{theorem}  If
$u:\RR^{n-1}\to \RR$ is an entire solution to the minimal surface
equation and $n\leq 8$, then $u$ is an affine function.
\end{theorem}

However, in 1969 E. Bombieri, De Giorgi, and E. Giusti
constructed entire non--affine solutions to the minimal surface
equation on $\RR^8$ and an area--minimizing singular cone in
$\RR^8$. In fact, they showed that for $m \geq 4$ the cones
\begin{equation}        \label{e:cones1}
        C_m = \{ (x_1 , \dots , x_{2m} ) \mid  x_1^2 + \cdots + x_m^2 =
x_{m+1}^2 + \cdots + x_{2m}^2  \} \subset \RR^{2m} \,
\end{equation}
are area--minimizing (and obviously singular at the origin).

\vskip2mm
 In contrast to the entire case, exterior solutions of the minimal graph
 equation, i.e., solutions on $\RR^2 \setminus B_1$, are much more
 plentiful. In this case,   L. Bers proved that $\nabla u$
actually has an asymptotic limit:

\begin{theorem} \label{t:bers}
If $u$ is a $C^2$ solution to the minimal surface equation on
$\RR^2 \setminus B_1$, then $\nabla u$ has a limit at infinity
(i.e., there is an asymptotic tangent plane).
\end{theorem}

 Bers' theorem was extended to higher
dimensions by L. Simon:

\begin{theorem}   \label{t:berssi}
If $u$ is a $C^2$ solution to the minimal surface equation on
$\RR^n \setminus B_1$, then either \begin{itemize} \item $|\nabla
u|$ is bounded and $\nabla u$ has a limit at infinity. \item All
tangent cones at infinity are of the form $\Sigma \times \RR$
where $\Sigma$ is singular. \end{itemize}
\end{theorem}

Bernstein's theorem has had many other interesting
generalizations, some of which will be discussed later.

\section{Simons inequality}  \label{s:9}

In this section, we recall a very useful differential inequality
for the Laplacian of the norm squared of the second fundamental
form of a minimal hypersurface $\Sigma$ in $\RR^n$ and illustrate
its role in a priori estimates. This inequality, originally due to
J. Simons, is:

\begin{lemma}      \label{l:simonsine}
If $\Sigma^{n-1}\subset \RR^n$ is a minimal hypersurface, then
\begin{equation}    \label{e:simtype}
\Delta_{\Sigma} \, |A|^2  = - 2\, |A|^4 + 2 |\nabla_{\Sigma} A|^2
\geq - 2 \, |A|^4 \, .
\end{equation}
\end{lemma}

An inequality of the type \eqr{e:simtype} on its own does not lead
to pointwise bounds on $|A|^2$ because of the nonlinearity.
However, it does lead to estimates if a   ``scale--invariant
energy'' is small.  For example, H. Choi and Schoen used
\eqr{e:simtype}  to prove:

\begin{theorem}       \label{t:cisc}
There exists $\epsilon > 0$ so that if $0 \in \Sigma \subset B_r
(0)$ with $\partial \Sigma \subset
\partial B_r (0)$ is a minimal surface with
\begin{equation}
    \int |A|^2 \leq \epsilon \, ,
\end{equation}
 then
\begin{equation}
    |A|^2(0) \leq   r^{-2} \, .
\end{equation}
\end{theorem}

\section{Heinz's curvature estimate for graphs}

One of the key themes in minimal surface theory is the usefulness
of a priori estimates.
   A basic example is the curvature estimate of
E. Heinz for graphs. Heinz's estimate gives an effective version
of the Bernstein's theorem; namely, letting the radius $r_0$ go to
infinity in \eqr{e:heinz} implies that $|A|$ vanishes, thus giving
Bernstein's theorem.

\begin{theorem}            \label{t:heinz}
If $D_{r_0} \subset \RR^2$ and
 $u:D_{r_0}\to \RR$ satisfies
the minimal surface equation, then for $\Sigma={\text{Graph}}_u$
and $0<\sigma\leq r_0$
\begin{equation}        \label{e:heinz}
        \sigma^2\,\sup_{D_{r_0-\sigma}}|A|^2\leq C\, .
\end{equation}
\end{theorem}

\begin{proof}
(Sketch.)  Observe first that it suffices to prove the estimate
for $\sigma = r_0$, i.e., to show that
\begin{equation}        \label{e:heinz2}
        |A|^2 (0,u(0)) \leq C \, r_0^{-2} \, .
\end{equation}
Recall that minimal graphs are automatically stable. As in the
proof of Theorem \ref{t:bern}, the area estimate for graphs
\eqr{e:area} allows us to use a logarithmic cutoff function in the
the stability inequality \eqr{e:stabin} to get that
\begin{equation}
    \int_{B_{r_1} \cap \Graph_u}  |A|^2 \leq \frac{ C }{ \log (r_0/r_1) } \, .
\end{equation}
Taking $r_0/ r_1$ sufficiently large, we can then apply Theorem
\ref{t:cisc} to get \eqr{e:heinz2}.
\end{proof}

\section{Embedded minimal disks with area bounds}

In the early nineteen--eighties Schoen and Simon extended the
theorem of Bernstein to complete simply connected embedded minimal
surfaces in $\RR^3$ with quadratic area growth. A
surface $\Sigma$ is said to have quadratic area growth if for all
$r>0$, the intersection of the surface with the ball in $\RR^3$ of
radius $r$ and center at the origin is bounded by $C\, r^2$ for a
fixed constant $C$ independent of $r$.

\begin{theorem}      \label{t:quasgm}
Let $0 \in \Sigma^{2} \subset B_{r_0} = B_{r_0}(x) \subset \RR^3$
be an embedded simply connected minimal surface with $\partial
\Sigma \subset \partial B_{r_0}$. If $\mu > 0$ and either
\begin{equation}
        \Area (\Sigma) \leq \mu \, r_0^2 \,  {\text{ or }} \,
        \int_{\Sigma} |A|^2 \leq \mu \, ,
\end{equation}
then for the connected component $\Sigma'$ of $B_{r_0 / 2}(x_0)
\cap \Sigma$ with $0 \in \Sigma'$ we have
\begin{equation}        \label{e:oqgmc}
        \sup_{\Sigma'}  |A|^2
                \leq C \, r_0^{-2}
\end{equation}
for some $C=C(\mu)$.
\end{theorem}

 The result of Schoen--Simon was generalized by Colding--Minicozzi
to quadratic area growth for
{\it{intrinsic}} balls (this generalization played an important
role in analyzing the local structure of embedded minimal
surfaces):

\begin{theorem}
Given a constant $C_I$, there exists $C_P$ so that if
$\cB_{2r_0}\subset \Sigma\subset \RR^3$ is an embedded minimal
disk satisfying either
\begin{equation}
        \Area (\cB_{2r_0}) \leq C_I \, r_0^2 \,  {\text{ or }}
\,
        \int_{\cB_{2r_0}} |A|^2 \leq C_I \, ,
\end{equation}
then
\begin{equation}
    \sup_{\cB_{s}} |A|^2 \leq C_P \, s^{-2} \,
.
\end{equation}
\end{theorem}

As an immediate consequence, letting $r_0 \to \infty$ gives
Bernstein-type theorems for embedded simply connected minimal
surfaces with either bounded density or finite total curvature.
Note that Enneper's surface is simply-connected but neither flat
nor embedded; this shows that embeddedness is essential for these
estimates. Similarly, the catenoid  shows that simply-connected is
essential. The catenoid is the minimal surface in $\RR^3$ given by
\begin{equation}
    \{ (\cosh s\, \cos t,\cosh s\, \sin t,s) \, | \,
    s,t\in\RR \} \, .
\end{equation}

\section{Stable minimal surfaces}

It turns out that stable minimal surfaces have a priori estimates.
 Since minimal graphs are stable,
the estimates for stable surfaces can be thought of as
generalizations of the earlier estimates for graphs.  These
estimates have been widely applied and are particularly useful
when combined with existence results for stable surfaces (such as
the solution of the Plateau problem).
  The starting point for these estimates is that, as we saw in
\eqr{e:secvar}, stable
minimal surfaces satisfy
 the
stability inequality
\begin{equation}   \label{e:stabin2}
      \int |A|^2 \, \phi^2 \leq   \int |\nabla \phi|^2 \, .
\end{equation}

We will mention two such estimates.  The first is R. Schoen's
curvature estimate for stable surfaces:

\begin{theorem}   \label{t:stable2}
There exists a constant $C$ so that if $\Sigma  \subset  \RR^3$ is
an immersed stable minimal surface with trivial normal bundle and
$\cB_{r_0} \subset \Sigma \setminus \partial \Sigma$, then
\begin{equation}        \label{e:o11}
        \sup_{ \cB_{r_0-\sigma} }  |A|^2
                \leq C \, \sigma^{-2}\, .
\end{equation}
\end{theorem}

The second is an estimate for the area and total curvature of a
stable surface is due to Colding--Minicozzi; for simplicity,
we will state only the area
estimate:

\begin{theorem}    \label{t:stable3}
If $\Sigma  \subset  \RR^3$ is an immersed stable minimal surface
with trivial normal bundle and $\cB_{r_0} \subset \Sigma \setminus
\partial \Sigma$, then
\begin{equation}        \label{e:abb}
       \Area \, (\cB_{r_0}) \leq 4 \pi \, r_0^2 / 3 \, .
\end{equation}
\end{theorem}

As mentioned,  we can use \eqr{e:abb} to bound the energy of a
 cutoff function in the stability inequality and, thus,
bound the total curvature of sub-balls.  Combining this with the
curvature estimate of Theorem \ref{t:cisc} gives Theorem
\ref{t:stable2}.   Note that the bound \eqr{e:abb}
is surprisingly sharp; even when $\Sigma$ is a plane, the area is
$\pi r_0^2$.

\section{Regularity theory}        \label{s:18}

In this section, we survey some of the key  ideas in classical
regularity theory, such as the role of monotonicity, scaling,
$\epsilon$-regularity theorems (such as Allard's theorem) and
tangent cone analysis (such as Almgren's refinement of Federer's
dimension reducing).  We refer to the book \cite{Mo} for a more
detailed overview and a general introduction to geometric measure
theory.

The starting point for all of this is the monotonicity of volume
for a minimal $k$--dimensional submanifold $\Sigma$. Namely,
Corollary \ref{c:cmon}  gives that the density
\begin{equation}                \label{e:thetadefvv}
        \Theta_{x_0}(s)
        =\frac{\Vol (B_s (x_0) \cap \Sigma)}{\Vol (B_s \subset \RR^k)
        }
\end{equation}
is a monotone non-decreasing function of $s$.  Consequently, we
can define the density $\Theta_{x_0}$ at the point $x_0$ to be the
limit as $s \to 0$ of $\Theta_{x_0}(s)$.  It  also follows easily
from monotonicity that the density is semi--continuous as a
function of $x_0$.

\vskip2mm \subsection{$\epsilon$--regularity and the singular set}

An $\epsilon$--regularity theorem is a theorem  giving that a weak
(or generalized) solution is actually smooth at a point if a
scale--invariant energy is small enough there. The standard
example is the Allard regularity theorem:

\begin{theorem} \label{t:allard}
There exists $\delta (k , n) > 0$ such that if $\Sigma \subset
\RR^{n}$ is a $k$--rectifiable stationary varifold (with density
at least one a.e.), $x_0 \in \Sigma$, and
\begin{equation}
    \Theta_{x_0} = \lim_{r\to 0} \frac{ \Vol \, (B_r (x_0) \cap \Sigma)}
        {\Vol \, ( B_r \subset \RR^k )}
        < 1 + \delta \, ,
\end{equation}
 then $\Sigma$ is smooth in a neighborhood of $x_0$.
\end{theorem}

Similarly, the small total curvature estimate of Theorem
\ref{t:cisc} may be thought of as an $\epsilon$-regularity
theorem; in this case, the scale--invariant energy is $\int
|A|^2$.

\vskip2mm
 As an application of the  $\epsilon$--regularity
theorem, Theorem \ref{t:allard},  we can define the singular set
$\cS$ of $\Sigma$ by
\begin{equation}
    \cS = \{ x \in \Sigma \, | \, \Theta_x \geq 1 + \delta \} \, .
\end{equation}
It follows immediately from the semi--continuity of the density
  that $\cS$ is closed.  In order to bound the
size of the singular set (e.g., the Hausdorff measure), one
combines the $\epsilon$--regularity with simple covering
arguments.

  This preliminary analysis of the singular set can be
refined by doing a so--called tangent cone analysis.

\subsection{Tangent cone analysis}

It is not hard to see that scaling  preserves the space of
  minimal
submanifolds of $\RR^n$.  Namely, if $\Sigma$ is minimal, then so
is
\begin{equation}
     \Sigma_{y,\lambda}  = \{ y + \lambda^{-1}\, (x-y) \, | \, x \in \Sigma \} \, .
\end{equation}
(To see this, simply note that this scaling multiplies the
principal curvatures by $\lambda$.) Suppose now that we fix the
point $y$ and take a sequence $\lambda_j \to 0$.  The monotonicity
formula   bounds the density of the rescaled solution, allowing us
to extract a convergent subsequence and limit.  This limit, which
is called a {\it tangent cone} at $y$,   achieves equality in the
monotonicity formula and, hence, must be homogeneous (i.e.,
invariant under dilations about $y$).

The usefulness of tangent cone analysis in
 regularity theory is based on two key facts.  For simplicity, we
 illustrate these when $\Sigma \subset \RR^n$ is an area
 minimizing hypersurface.   First, if any tangent cone at $y$ is a
 hyperplane $\RR^{n-1}$, then $\Sigma$ is smooth in a neighborhood
 of $y$.  This follows easily from the Allard regularity theorem since
 the density at $y$ of the tangent cone is the same as the density
 at $y$ of $\Sigma$.  The second key fact, known as ``dimension
 reducing,'' is due to Almgren
  and is a refinement of an argument of
 Federer.  To state this, we first stratify the singular set $\cS$
 of $\Sigma$ into subsets
 \begin{equation}
    \cS_0 \subset \cS_1 \subset \cdots \subset \cS_{n-2} \, ,
 \end{equation}
 where we define $\cS_i$ to be the set of points $y\in \cS$ so
 that  any linear space contained in any tangent cone at $y$ has
 dimension at most $i$.  (Note that $\cS_{n-1} = \emptyset$ by Allard's
 Theorem.)  The dimension reducing argument then gives that
 \begin{equation}       \label{e:dimred}
    {\text{dim}} \, (\cS_i) \leq i \, ,
 \end{equation}
 where dimension   means the Hausdorff dimension.  In particular, the
 solution of the Bernstein problem then gives codimension $7$
 regularity of $\Sigma$, i.e., ${\text{dim}} \, (\cS) \leq n-8$.

\part{Constructing minimal surfaces}   \label{p:4}

Thus far, we have mainly dealt with regularity and a priori
estimates but have ignored questions of existence.  In this part
we survey some of the most useful existence results for minimal
surfaces. Section \ref{s:plateau}   gives an overview of the
classical Plateau problem.  Section \ref{s:21} recalls the
classical Weierstrass representation, including a few modern
applications, and the Kapouleas desingularization method.  Section
\ref{s:22} deals with producing area minimizing surfaces  and
questions of embeddedness. Finally, Section \ref{s:24} recalls the
min--max construction for producing unstable minimal surfaces and,
in particular, doing so while controlling the topology and
guaranteeing embeddedness.

\section{The Plateau Problem}  \label{s:plateau}

The following fundamental existence problem for minimal surfaces
is known as the Plateau problem:   Given a closed curve $\Gamma$,
find a minimal surface with boundary $\Gamma$. There are various
solutions to this problem depending on the exact definition of a
surface (parameterized disk, integral current, $\ZZ_2$ current, or
rectifiable varifold).
 We shall consider the
version of the Plateau problem for parameterized disks; this was
solved independently by J. Douglas   and T. Rado. The
generalization to Riemannian manifolds is due to C. B. Morrey.

\begin{theorem}     \label{t:plat}
Let $\Gamma \subset \RR^3$ be a piecewise $C^1$ closed Jordan
curve. Then there exists a piecewise $C^1$ map $u$ from $D \subset
\RR^2$ to $\RR^3$ with $u(\partial D) \subset \Gamma$ such that
the image  minimizes area among all disks with boundary $\Gamma$.
\end{theorem}

The solution $u$ to the Plateau problem  above can easily be seen
to be a branched conformal immersion.
 R. Osserman proved that $u$ does not have true interior branch
 points; subsequently,
R. Gulliver and W. Alt showed that $u$ cannot have false branch
points either.

Furthermore, the solution $u$ is as smooth as the boundary curve,
even up to the boundary.    A very general version of this
boundary regularity was proven by S. Hildebrandt; for the case of
surfaces in $\RR^3$, recall the following result of J. C. C.
Nitsche:

\begin{theorem}              \label{t:nitsche2}
If $\Gamma$ is a regular Jordan curve of class $C^{k,\alpha}$
where $k \geq 1$ and $0 < \alpha < 1$, then a solution $u$ of the
Plateau problem is $C^{k,\alpha}$ on all of $\bar{D}$.
\end{theorem}

\section{The Weierstrass representation}        \label{s:21}

The classical Weierstrass representation (see \cite{Os}) takes
holomorphic data (a Riemann surface, a meromorphic function, and a
holomorphic one--form) and associates a minimal surface in
$\RR^3$.
 To be precise, given
a Riemann surface $\Omega$, a meromorphic function $g$ on
$\Omega$, and a holomorphic one--form $\phi$ on $\Omega$, then we
get  a (branched) conformal minimal immersion $F: \Omega \to
\RR^3$ by
\begin{equation}    \label{e:ws1}
    F(z) = {\text{Re }} \int_{\zeta \in \gamma_{z_0,z}}
\left( \frac{1}{2} \, (g^{-1} (\zeta) - g (\zeta) )
    , \frac{i}{2} \, (g^{-1} (\zeta)
    +g (\zeta) ) , 1 \right) \, \phi (\zeta) \, .
\end{equation}
Here $z_0 \in \Omega$ is a fixed base point and the integration is
along a path $\gamma_{z_0,z}$ from $z_0$ to $z$. The choice of
$z_0$ changes $F$ by adding a constant. In general, the map $F$
may depend on the choice of path (and hence may not be
well--defined);  this is known as ``the period problem''. However,
when $g$ has no zeros or poles and $\Omega$ is simply connected,
then
 $F(z)$ does not depend on the choice of path
$\gamma_{z_0,z}$.

Two standard constructions of minimal surfaces from Weierstrass
data are
 \begin{align}
 &g (z) = z, \, \phi (z)  = dz/z ,
\, \Omega = \CC \setminus \{ 0 \} {\text{ giving a catenoid}}
\, , \\
  &g (z) = \e^{iz} , \,
\phi (z) = dz , \, \Omega = \CC  {\text{ giving a helicoid}} \, .
\label{e:hel}
\end{align}

The Weierstrass representation is particularly useful for
constructing immersed minimal surfaces.
 Typically, it is rather difficult to prove that the resulting
immersion is an embedding (i.e., is $1$--$1$), although there are
some interesting cases where this can be done.    For the first
modern example, D. Hoffman and Meeks proved that the surface
constructed by Costa was embedded; this was the first new complete
finite topology properly embedded minimal surface discovered since
the classical catenoid, helicoid, and plane. This led to the
discovery of many more such surfaces
 (see \cite{Ro} for more discussion).

\section{Area--minimizing surfaces}     \label{s:22}

Perhaps the most natural way to construct minimal surfaces is to
look for ones which minimize area, e.g., with fixed boundary, or
in a homotopy class, etc.  This has the advantage that often it is
possible to show that the resulting surface is embedded.  We
mention a few results along these lines.

The first embeddedness result, due to   Meeks and Yau, shows that
if the boundary curve is embedded and lies on the boundary of a
smooth mean convex set (and it is null--homotopic in this  set),
then it bounds an embedded least area disk.

\begin{theorem}  \cite{MeYa1}             \label{t:my1}
Let $M^3$ be a compact Riemannian three--manifold whose boundary
is mean convex and let $\gamma$ be a simple closed curve in
$\partial M$ which is null--homotopic in $M$; then $\gamma$ is
bounded by a least area disk and any such least area disk is
properly embedded.
\end{theorem}

Note that some restriction on the boundary curve $\gamma$ is
certainly necessary. For instance, if the boundary curve was
knotted (e.g., the trefoil), then it could not be spanned by any
embedded disk (minimal or otherwise).  Prior to the work of Meeks
and Yau, embeddedness was known for extremal boundary curves in
$\RR^3$ with small total curvature by the work of R. Gulliver and
J. Spruck.

 If we instead fix a homotopy class of maps, then the two fundamental
existence results
 are due to Sacks--Uhlenbeck and Schoen--Yau (with embeddedness proven
 by Meeks--Yau and Freedman--Hass--Scott,
 respectively):

\begin{theorem}     \label{p:existence1}
Given $M^3$, there exist conformal  (stable) minimal immersions
$u_1 , \dots , u_m  : \SS^2 \to M$ which generate $\pi_2 (M)$ as a
$\ZZ[\pi_1 (M)]$ module.  Furthermore,
\begin{itemize}
\item
If $u: \SS^2 \to M$ and $[u]_{\pi_2} \ne 0$, then $\Area (u) \geq
\min_i \Area (u_i)$.
 \item Each $u_i$ is either an
embedding or a $2$--$1$ map onto an embedded $2$--sided $\RP^2$.
\end{itemize}
\end{theorem}

\begin{theorem}     \label{p:existence2}
If $\Sigma^2$ is a closed surface with genus $g>0$ and $i_0 :
\Sigma \to M^3$ is an embedding which induces an injective map on
$\pi_1$, then there is a least area embedding with the same action
on $\pi_1$.
\end{theorem}

\section{The min--max construction of minimal surfaces} \label{s:24}

Variational arguments can also be used to construct higher index
(i.e., non--minimizing) minimal surfaces using   the topology of
the space of surfaces. There are two basic approaches:
\begin{itemize}
\item
Applying Morse theory to the energy functional on the space of
maps from a fixed surface $\Sigma$ to $M$.
\item
Doing a min--max argument over  families of (topologically
non--trivial) sweep--outs of $M$.
\end{itemize}
The first approach has the advantage that the topological type of
the minimal surface is easily fixed; however, the second approach
has been  more successful at producing embedded minimal surfaces.
We will highlight a few key  results below but refer to \cite{CD}
for a thorough treatment.

Unfortunately, one cannot directly apply Morse theory to the
energy functional on the space of maps from a fixed surface
because of a lack of compactness (the Palais--Smale Condition C
does not hold).  To get around this difficulty, Sacks-Uhlenbeck
introduce a family of perturbed energy functionals which do
satisfy Condition C and then obtain   minimal surfaces as limits
of critical points for the perturbed problems:

\begin{theorem}      \label{t:nonasph}
If   $\pi_k (M) \ne 0$ for some $k>1$, then there exists a
branched immersed minimal $2$--sphere in $M$ (for any metric).
\end{theorem}

\begin{figure}[htbp]
    \setlength{\captionindent}{4pt}
    \begin{minipage}[t]{0.75\textwidth}
    \centering\input{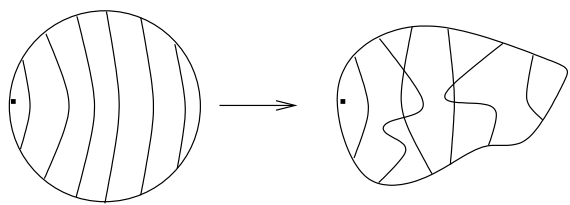}
    \caption{A $1$--parameter family of curves on a $2$--sphere which
induces a map $F:\SS^2 \to
    \SS^2$ of degree $1$. First published in Surveys in Differential Geometry,
    volume IX, in 2004, published by International Press.}  \label{f:fsweep}
    \end{minipage}
\end{figure}

\vskip2mm The basic idea of constructing minimal surfaces via
min--max arguments and sweep--outs goes back to Birkhoff, who
developed it to construct simple closed geodesics on spheres.  In
particular, when $M$ is a topological $2$--sphere, we can find a
$1$--parameter family of curves starting and ending at point
curves so that the induced map $F:\SS^2 \to \SS^2$ (see fig.
\ref{f:fsweep}) has nonzero degree.    The min--max argument
produces a nontrivial closed geodesic of length less than or equal
to the longest curve in the initial one--parameter family.  A
curve shortening argument gives that the geodesic obtained in this
way is simple.

J. Pitts applied a similar argument and   geometric
measure theory to get that every closed Riemannian three manifold
has an embedded minimal surface (his argument was for dimensions
up to seven), but he did not estimate the genus of the resulting
surface.  Finally, F. Smith (under the direction of L. Simon)
proved (see \cite{CD}):

\begin{theorem}    \label{t:smith}
Every metric   on  a topological $3$--sphere $M$ admits an
embedded minimal $2$--sphere.
\end{theorem}

The main new contribution of Smith was to control the topological
type of the resulting minimal surface while keeping it embedded.

\part{Some applications of minimal surfaces}

\stepcounter{section}

In this last part, we discuss very briefly a few applications of
minimal surfaces.  As mentioned in the introduction, there are
many to choose from and we have selected just a few.

\subsection{The positive mass theorem}  \label{s:posmass}

The (Riemannian version of the) positive mass theorem states that
an asymptotically flat $3$-manifold $M$ with non-negative scalar
curvature must have positive mass.  The Riemannian manifold $M$
here arises as a maximal space-like slice in a $3+1$-dimensional
space-time solution of Einstein's equations.

The asymptotic flatness of $M$ comes from that the space-time
models an isolated gravitational system and hence is a
perturbation of the vacuum solution outside a large compact set.
 To make this precise, suppose for simplicity that $M$ has only
 one end;
 $M$ is then said to be {\it asymptotically flat} if there is a compact
set $\Omega \subset M$ so that $M \setminus \Omega$ is
diffeomorphic to $\RR^3 \setminus B_R (0)$ and the metric on $M
\setminus \Omega$ can be written as
\begin{equation}    \label{e:pdef}
    g_{ij} = \left(1 + \frac{\cM}{ 2 \, |x|} \right)^4 \, \delta_{ij} + p_{ij} \, ,
\end{equation}
where
\begin{equation}    \label{e:pbds}
    |x|^2 \, |p_{ij}| + |x|^3 \, |D \, p_{ij}|
    + |x|^4 \, |D^2 \, p_{ij}|  \leq C \, .
\end{equation}
The constant $\cM$ is the so called {\it mass} of $M$.  Observe
that the metric $g_{ij}$ is a perturbation of the metric on a
constant-time slice in the Schwarzschild space-time of mass $\cM$;
that is to say, the Schwarzschild metric has $p_{ij} \equiv 0$.

A tensor $h$ is said to be $O(|x|^{-p})$ if $|x|^p \, |h| +
|x|^{p+1} \, |D \, h| \leq C$.  For example, an easy calculation
shows that
\begin{align}   \label{e:pdefst}
    g_{ij} &= \left(1 + 2 \, \cM / |x| \right) \, \delta_{ij} +O(|x|^{-2}) \, , \\
    \sqrt{g} &\equiv \sqrt{\det g_{ij} } = 1 + 3 \, \cM \, |x|^{-1} +
    O(|x|^{-2}) \, . \notag
\end{align}

The positive mass theorem states that the mass $\cM$ of such an
$M$ must be non-negative:

\begin{theorem} \label{t:pm1}
\cite{ScYa}  With $M$ as above, $\cM \geq 0$.
\end{theorem}

There is a rigidity theorem as well which states that the mass
vanishes only when $M$ is isometric to $\RR^3$:

\begin{theorem} \label{t:pm2}
\cite{ScYa} If $|\nabla^3 p_{ij}|  = O(|x|^{-5})$
 and  $\cM = 0$ in Theorem \ref{t:pm1},
then $M$ is {\underline{isometric}} to $\RR^3$.
\end{theorem}

 We will give a very brief overview of the proof of
 Theorem \ref{t:pm1}, showing in the process where minimal surfaces appear.

\begin{proof}
(Sketch.) The argument will by contradiction, so suppose that the
mass is negative.  It is not hard to prove that the slab between
two parallel planes is mean-convex.  That is, we have the
following:

\begin{lemma} \label{l:ecalc}
If $\cM < 0$ and $M$ is asymptotically flat, then there exist $R_0
, h > 0$ so that for $r > R_0$ the sets
\begin{equation}
    C_r = \{ |x|^2  \leq r^2  \, , \, -h \leq x_3 \leq h \}
\end{equation}
have strictly mean convex boundary.
\end{lemma}

Since the compact set $C_r$ is mean convex, we can solve the
Plateau problem (as in Section \ref{s:plateau}) to get an area
minimizing (and hence stable) surface $\Gamma_r \subset C_r$ with
boundary
\begin{equation}
\partial \Gamma_r = \{
|x|^2 = r^2 \, , \, x_3 = h \} \, .
\end{equation}
  Using the disk
$\{ |x|^2 \leq r^2 \, , x_3 = h \}$ as a comparison surface, we
get uniform local area bounds for any such $\Gamma_r$. Combining
these local area bounds with the a priori curvature estimates for
minimizing surfaces, we can take a sequence of $r$'s going to
infinity and find a subsequence of $\Gamma_r$'s that converge to a
{\underline{complete}} area-minimizing surface
\begin{equation}
    \Gamma \subset \{ - h \leq x_3 \leq h \} \, .
\end{equation}
Since $\Gamma$ is pinched between the planes $\{ x_3 = \pm h\}$,
the estimates for minimizing surfaces implies that (outside a
large compact set) $\Gamma$ is a graph over the plane $\{ x_3 = 0
\}$ and hence has quadratic area growth and finite total
curvature.  Moreover, using the form of the metric $g_{ij}$, we
see that $|\nabla u|$ decays like $|x|^{-1}$ and
\begin{equation}    \label{e:kg1}
    \int_{\sigma_s} k_g  = (2 \, \pi \, s + O(1)) \, (s^{-1} + O(s^{-2}))
    = 2 \, \pi + O(s^{-1}) \, ,
\end{equation}
where $\sigma_s = \{ x_1^2 + x_2^2 = s^2 \} \cap \Gamma$ and $k_g$
is the geodesic curvature of $\sigma_s$ (as a curve in $\Gamma$).

To get the contradiction, one combines stability of $\Gamma$ with
the positive scalar curvature of $M$ to see that no such $\Gamma$
could have existed.{\footnote{$M$ was assumed only to have
non-negative scalar curvature.  However, a ``rounding off''
argument shows that the metric on $M$ can be perturbed to have
positive scalar curvature outside of a compact set and still have
negative mass.}} Namely,  substituting the Gauss equation into the
stability inequality{\footnote{This is the stability inequality in
a general $3$-manifold; see \cite{CM1}.}}
 gives
\begin{equation}    \label{e:stabipm}
\int_{\Gamma} (|A|^2 / 2+  \Scal_M -   K_{\Sigma} ) \phi^2 \leq \,
\int_{\Gamma} |\nabla \phi|^2  \, .
\end{equation}
Since  $\Gamma$ has quadratic area growth, we can choose a
sequence of (logarithmic) cutoff functions in \eqr{e:stabipm} to
get
\begin{equation}    \label{e:stabipm2}
0 < \int_{\Sigma} (|A|^2 / 2 +  \Scal_M )
      \leq  \int_{\Sigma}   K_{\Sigma} < \infty \, ;
\end{equation}
since $K_{\Sigma}$ may not be positive, we also used that $\Gamma$
has finite total curvature.  Moreover, we used that $\Scal_M$ is
positive outside a compact set to see that the first integral  in
\eqr{e:stabipm2} was positive.   Finally,  substituting
\eqr{e:stabipm2} into the Gauss-Bonnet formula gives that
$\int_{\sigma_s} k_g$ is {\underline{strictly}} less than $2 \pi$
for $s$ large, contradicting \eqr{e:kg1}.
\end{proof}

\subsection{Black holes}

Another way that minimal surfaces enter into relativity is through
black holes.   Suppose that we have a three-dimensional time-slice
$M$ in a $3+1$-dimensional space--time. For simplicity, assume
that $M$ is totally geodesic and hence has non-negative scalar
curvature.
  A  closed surface $\Sigma$ in $M$ is  said
to be trapped if its mean curvature is everywhere negative with
respect to its outward normal.  Physically, this means that the
surface emits an outward shell of light whose surface area is
decreasing everywhere on the surface.  The existence of a closed
trapped surface implies the existence of a black hole in the
space-time.

 Given a
trapped surface, we can look for the outermost trapped surface
containing it; this outermost surface is called an apparent
horizon.  It is not hard to see that an apparent horizon must be a
minimal surface and, moreover, a barrier argument shows that it
must be stable. Since $M$ has non-negative scalar curvature,
stability in turn implies that it must be diffeomorphic to a
sphere.  See, for instance, \cite{Br} for references to some
results on black holes, horizons, etc.

\subsection{Constant mean curvature surfaces}

At least since the time of Plateau, minimal surfaces have been
used to model soap films.  This is because the mean curvature of
the surface models the surface tension and this is essentially the
only force acting on a soap film. Soap bubbles, on other hand,
enclose a volume and thus the pressure gives a second
counterbalancing force.  It follows easily that these two forces
are in equilibrium when the surface has constant mean curvature.

For the same reason,   constant mean curvature surfaces arise in
the isoperimetric problem.  Namely, a surface that minimizes
surface area while enclosing a fixed volume must have constant
mean curvature (or ``cmc'').  It is not hard to see that such an
isoperimetric surface in $\RR^n$ must be a round sphere.   There
are two interesting partial converses to this.  First, by a
theorem of Hopf, any cmc $2$-sphere in $\RR^3$   must be round.
Second, using the maximum principle (``the method of moving
planes'') Alexandrov showed that any closed embedded cmc
hypersurface in $\RR^n$ must be a round sphere. It turned out,
however, that not every closed immersed cmc surface is round.  The
first examples were immersed cmc tori constructed by H. Wente.
Kapouleas constructed many new examples, including closed higher
genus cmc surfaces.

Many of the techniques developed for studying minimal surfaces
generalize to general constant mean curvature surfaces.

\subsection{Finite extinction for Ricci flow}

We close this survey with indicating how minimal surfaces can be
used to show that on a homotopy $3$-sphere the Ricci flow become
extinct in finite time (see \cite{CM2}, \cite{Pe} for details).

Let $M^3$ be a smooth closed orientable $3$--manifold and let
$g(t)$ be a one--parameter family of metrics on $M$ evolving by
the Ricci flow, so
\begin{equation}  \label{e:eqRic}
 \partial_t g=-2\,\Ric_{M_t}\, .
\end{equation}

In an earlier section, we saw that there is a natural way of
constructing minimal surfaces on many $3$-manifolds and that comes
from the min--max argument where the minimal of all maximal slices
of sweep--outs is a minimal surface. The idea is then to look at
how the area of this min--max surface changes under the flow.
Geometrically the area measures a kind of width of the
$3$--manifold and as we will see for certain $3$--manifolds
(those, like the $3$--sphere, whose prime decomposition contains
no aspherical factors) the area becomes zero in finite time
corresponding to that the solution becomes extinct in finite time.

\begin{figure}[htbp]
\begin{center}
    \input{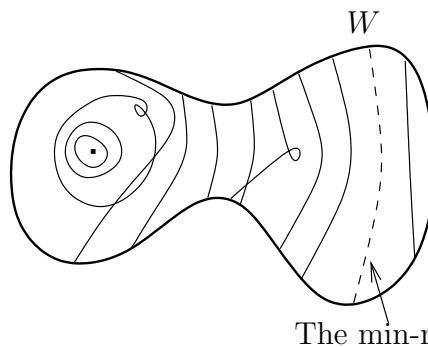}
    \caption{The sweep--out, the min--max surface, and the
    width W. First published in the Journal of the American Mathematical Society in 2005,
    published by the American Mathematical Society.}
    \label{f:1}
\end{center}
\end{figure}

Fix a continuous map $\beta: [0,1] \to C^0\cap L_1^2 (\SS^2 , M)$
where $\beta (0)$ and $\beta (1)$ are constant maps so that
$\beta$ is in  the nontrivial homotopy class $[\beta]$ (such
$\beta$ exists when $M$ is a homotopy $3$-sphere). We define the
width  $W=W(g,[\beta])$ by
\begin{equation}    \label{e:w3}
   W(g) = \min_{\gamma \in [ \beta]}
\, \max_{s \in [0,1]} \Energy (\gamma(s)) \,
   .
\end{equation}

\vskip2mm The next theorem gives an upper bound for the derivative
of $W(g(t))$ under the Ricci flow which forces the solution $g(t)$
to become extinct in finite.

\begin{theorem}     \label{t:upper}
Let $M^3$ be a homotopy $3$-sphere equipped with a Riemannian
metric $g=g(0)$. Under the Ricci flow, the width $W(g(t))$
satisfies
\begin{equation}   \label{e:di1a}
\frac{d}{dt} \, W(g(t))  \leq -4 \pi + \frac{3}{4 (t+C)} \,
W(g(t))   \, ,
\end{equation}
in the sense of the limsup of forward difference quotients. Hence,
 $g(t)$  must become extinct in finite time.
\end{theorem}

The $4\pi$ in \eqr{e:di1a} comes from the Gauss--Bonnet theorem
and the $3/4$ comes from the bound on the minimum of the scalar
curvature that the evolution equation implies.  Both of these
constants matter whereas the constant $C$ depends on the initial
metric and the actual value is not important.

To see that \eqr{e:di1a} implies finite extinction time rewrite
\eqr{e:di1a} as
\begin{equation}
\frac{d}{dt} \left( W(g(t)) \, (t+C)^{-3/4} \right) \leq - 4\pi \,
(t+C)^{-3/4}
\end{equation} and
integrate to get
\begin{equation}  \label{e:lastaa}
 (T+C)^{-3/4} \, W(g(T)) \leq C^{-3/4} \, W(g(0))
- 16 \, \pi \, \left[ (T+C)^{1/4} - C^{1/4} \right]   \, .
\end{equation}
Since $W \geq 0$ by definition and the right hand side of
\eqr{e:lastaa} would become negative for $T$ sufficiently large we
get the claim.

\vskip2mm As a corollary of this theorem we get finite extinction
time for the Ricci flow.

\begin{corollary}  \label{c:upper}
Let $M^3$ be a homotopy $3$-sphere equipped with a Riemannian
metric $g=g(0)$. Under the Ricci flow $g(t)$ must become extinct
in finite time.
\end{corollary}

\bibliographystyle{plain}

\begin{thebibliography}{A}







\bibitem[Br]{Br}
H. Bray, {\it Black holes, geometric flows, and the Penrose inequality in
general relativity}, Notices Amer. Math. Soc.  49  (2002),  no. 11, 1372--1381.




\bibitem[CD]{CD}
T.H. Colding and C. De Lellis, {\it The min--max construction of
minimal surfaces}, Surveys in differential geometry, Vol. 8,
Lectures on Geometry and Topology held in honor of Calabi, Lawson,
Siu, and Uhlenbeck at Harvard University, May 3--5, 2002,
Sponsored by the Journal of Differential Geometry, (2003) 75--107,
math.AP/0303305.



\bibitem[CM1]{CM1}
T.H. Colding and W.P. Minicozzi II, \textit{Minimal surfaces,
Courant Lecture Notes in Math.}, v. 4, 1999.
\bibitem[CM2]{CM2}
T.H. Colding and W.P. Minicozzi II, \textit{Estimates for the
extinction time for the Ricci flow on certain $3$--manifolds and a
question of Perelman}, JAMS,  18  (2005),  no. 3, 561--569,
math.AP/0308090.
\bibitem[CM3]{CM3}
T.H. Colding and W.P. Minicozzi II, \textit{Minimal submanifolds},
preprint.


\bibitem[CM4]{CM4}
T.H. Colding and W.P. Minicozzi II, \textit{Disks that are double
spiral staircases}, Notices of the AMS, Vol. 50, no. 3, March
(2003) 327--339.




























\bibitem[La]{La}
H. B. Lawson,  {\it Lectures on minimal submanifolds}, vol. I,
Publish or Perish, Inc, Berkeley, 1980.






\bibitem[MeP]{MeP}
W. Meeks III and J. Perez,   {\it Conformal properties in
classical minimal surface theory}, Surveys in Diff. Geom. IX:
Eigenvalues of Laplacians and other geometric operators, Ed. by A.
Grigor'yan and S.T. Yau, International Press (2004), 275--335.

\bibitem[MeYa1]{MeYa1}
W. Meeks  III and S.T. Yau,   {\it The classical Plateau problem
and the topology of three dimensional manifolds},
 Topology 21 (1982) 409--442.
\bibitem[MeYa2]{MeYa2}
W.H. Meeks and S.T. Yau,  \textit{Topology of three--dimensional
manifolds and the embedding problems in minimal surface theory},
Ann. of Math. (2) 112 (1980), no. 3, 441--484.


\bibitem[Mo]{Mo}
F. Morgan, {\it Geometric measure theory. A beginner's guide}.
Second edition. Academic Press, Inc., San Diego, CA, 1995.


\bibitem[Os]{Os}
R. Osserman, {\it A survey of minimal surfaces},   Dover, 2nd.
edition (1986).


\bibitem[Pe]{Pe}
G. Perelman, {\it Finite extinction time for the solutions to the
Ricci flow on certain three--manifolds},  math.DG/0307245.

\bibitem[Pz]{Pz}
J. Perez, {\it Limits by rescalings of minimal surfaces: Minimal
laminations, curvature decay and local pictures, notes for  the
workshop "Moduli Spaces of Properly Embedded Minimal Surfaces"},
American Institute of Mathematics, Palo Alto, California (2005).




\bibitem[Ro]{Ro}
H. Rosenberg, {\it Some recent developments in the theory of
properly embedded minimal surfaces in $\RR^3$}, Seminare Bourbaki
1991/92, Asterisque No. 206 (1992) 463--535.




 \bibitem[ScYa]{ScYa}
Schoen, R.  and Yau, S. T., {\it On the proof of the positive mass
conjecture in general relativity}, Comm. Math. Phys. 65 (1979),
no. 1, 45--76.






















\end{thebibliography}

\end{document}